\documentclass[12pt]{amsart}%
\usepackage{amsmath}
\usepackage{graphicx}
\usepackage{amscd}
\usepackage{geometry}
\usepackage{amsfonts}
\usepackage{amssymb}%
\newtheorem{theorem}{Theorem}[section]
\theoremstyle{plain}

\newtheorem{claim}{Claim}

\newtheorem{corollary}[theorem]{Corollary}

\newtheorem{example}{Example}

\newtheorem{lemma}[theorem]{Lemma}

\newtheorem{proposition}[theorem]{Proposition}
\newtheorem{remark}{Remark}

\numberwithin{equation}{section}

\begin{document}
\title[Products of open manifolds with $\mathbb{R}$ ]{Products of open manifolds with $\mathbb{R}$}
\author{Craig R. Guilbault }
\address{Department of Mathematical Sciences, University of Wisconsin-Milwaukee,
Milwaukee, Wisconsin 53201}
\email{craigg@uwm.edu}
\date{December 20, 2005}

\begin{abstract}
In this note we present a characterization of those open $n$-manifolds
($n\geq5$), whose products with the real line are homeomorphic to interiors of
compact $\left(  n+1\right)  $-manifolds with boundary.

\end{abstract}
\maketitle

\section{Introduction}

One often wishes to know whether a given open manifold can be compactified by
the addition of a manifold boundary. In other words, for an open manifold
$M^{n}$, we ask if there exists a compact manifold $C^{n}$ with $int\left(
C^{n}\right)  \approx M^{n}$. Since $int\left(  C^{n}\right)  \hookrightarrow
C^{n}$ is a homotopy equivalence, and because every compact manifold has the
homotopy type of a finite CW complex (see \cite{KS}), a necessary condition is
that $M^{n}$ have finite homotopy type. This condition is not sufficient. One
of the most striking illustrations of that fact occurs in a famous
contractible (thus, homotopy equivalent to a point) $3$-manifold constructed
by J.H.C. Whitehead \cite{Wh}. That example is best known for not being
homeomorphic $\mathbb{R}^{3}$ (or, equivalently, to $int\left(  B^{3}\right)
$), but a little additional thought reveals that it is not homeomorphic to the
interior of \emph{any} compact $3$-manifold.

Somewhat surprisingly, the product of the Whitehead manifold with a line is
homeomorphic to $\mathbb{R}^{4}$. In fact, it is now known that the product of
\emph{any} contractible $n$-manifold with a line is homeomorphic to
$\mathbb{R}^{n+1}$. That fact was obtained through the combined efforts of
several researchers; see, for example, \cite{Gl}, \cite{Mc}, \cite{St},
\cite{Lu1}, \cite{Lu2} and \cite{Fr} . In this note we prove the following
generalization of that result:

\begin{theorem}
\label{main theorem}For an open manifold $M^{n}$ ($n\geq5$), $M^{n}%
\times\mathbb{R}$ is homeomorphic to the interior of a compact $\left(
n+1\right)  $-manifold with boundary if and only if $M^{n}$ has the homotopy
type of a finite complex.\medskip
\end{theorem}

I wish to acknowledge Igor Belegradek for motivating this work by asking me
the question:

\begin{quotation}
\emph{If }$M^{n}$\emph{ is an open manifold homotopy equivalent to an embedded
compact submanifold, say a torus, must }$M^{n}\times R$\emph{ be homeomorphic
to the interior of a compact manifold?}
\end{quotation}

\noindent Initially, I was surprised that the question was open. The fairly
obvious approach---application of the main result of Siebenmann's
thesis---works nicely for $M^{n}\times\mathbb{R}^{2}$. In fact, Siebenmann,
himself, addressed that situation in his thesis \cite[Th.6.12]{Si}; where a
key ingredient was the straightforward observation that (for any connected
open manifold $M^{n}$), $M^{n}\times\mathbb{R}^{2}$ has stable fundamental
group at infinity. We too obtain our result by applying Siebenmann's thesis;
but unlike the `cross $\mathbb{R}^{2}$ situation', stability of the
fundamental group at infinity for $M^{n}\times\mathbb{R}$ is not so easy. In
fact, there exist open manifolds $M^{n}$ for which $M^{n}\times\mathbb{R}$
fails to have stable fundamental group at infinity. An example of that
phenomenon will be provided in Section \ref{Section: Proof}. However, under
the (already necessary) hypothesis of finite homotopy type---or even a weaker
hypothesis of finite domination---we are able to obtain $\pi_{1}$-stability in
$M^{n}\times\mathbb{R}$. That is the main step in our proof. A key ingredient
is the adaptation of a recent technique from \cite{GuTi}.

\section{Definitions and Background}

Throughout this paper, we work in the PL\ category. Proofs can be modified in
the usual ways to obtain equivalent results in the smooth or topological categories.

\subsection{Neighborhoods of infinity, ends, and finite dominations}

A manifold $M^{n}$ is \emph{open} if it is noncompact and has no boundary. A
subset $N$ of $M^{n}$ is a \emph{neighborhood of infinity} if $\overline
{M^{n}-N}$ is compact. We say that $M^{n}$ is \emph{one-ended} if each
neighborhood of infinity contains a connected neighborhood of infinity; in
other words, $M^{n}$ contains `arbitrarily small' connected neighborhoods of
infinity. More generally, $M^{n}$ is $k$\emph{-ended} ($k\in\mathbb{N}$) if it
contains arbitrarily small neighborhoods of infinity consisting of exactly $k$
components, each of which has noncompact closure. If no such $k$ exists, we
say $M^{n}$ has \emph{infinitely many ends}.

A neighborhood of infinity is \emph{clean} if it is a closed subset of $M^{n}$
and a codimension $0$ submanifold with a boundary that is bicollared in
$M^{n}$. By discarding compact components and drilling our neighborhoods arcs,
we can find within any clean neighborhood of infinity $N$, an \emph{improved
}clean neighborhood of infinity $N^{\prime}$ having the properties:

\begin{itemize}
\item $N^{\prime}$ contains no compact components, and

\item each component of $N^{\prime}$ has connected boundary.
\end{itemize}

\noindent If $M^{n}$ is $k$-ended then, there exist arbitrarily small improved
neighborhood of infinity containing exactly $k$ components. Such a
neighborhood is called $0$\emph{-neighborhood of infinity}. In this situation,
we may choose a sequence
\[
N_{1}\supseteq N_{2}\supseteq N_{3}\supseteq\cdots
\]
of $0$-neighborhoods of infinity such that $N_{i+1}\subseteq int\left(
N_{i}\right)  $ for all $i$ and $\cap_{i=1}^{\infty}N_{i}=\varnothing$. A
sequence of this sort will be referred to as \emph{neat}. Then for each $i$,
the components may be indexed as $N_{i}^{1},N_{i}^{2},\cdots,N_{i}^{k}$;
furthermore, these indices may chosen coherently so that for all $i<i^{\prime
}$ and all $1\leq j\leq k$ we have $N_{i^{\prime}}^{j}\subseteq int(N_{i}%
^{j})$. When all of the above has been accomplished, we will refer to
$\left\{  N_{i}\right\}  _{i=1}^{\infty}$ as a \emph{well-indexed neat
sequence of }$0$\emph{-neighborhoods of infinity. }For a fixed $j$, we say
that the nested sequence of components $\left\{  N_{i}^{j}\right\}
_{i=1}^{\infty}$ \emph{represents the }$j^{\emph{th}}$ \emph{end of }$M^{n}$.

A space has \emph{finite homotopy type }if it is homotopy equivalent to a
finite CW complex. A space $X$ is \emph{finitely dominated} if there exists a
finite complex $L$ and maps $u:X\rightarrow L$ and $d:L\rightarrow X$ such
that $d\circ u\simeq id_{X}$. It is a standard fact that a polyhedron (or
complex) $X$ is finitely dominated if and only if there is a homotopy
$H:X\times\left[  0,1\right]  \rightarrow X$ such that $H_{0}=id_{X}$ and
$\overline{H_{1}(X)}$ is compact. (We say $H$ \emph{pulls }$X$ \emph{into a
compact set}.) For later use, we prove a mild refinement of this latter characterization.

\begin{lemma}
\label{improved-domination}A polyhedron $X$ is finitely dominated if and only
if, for any compactum $C\subseteq X$, there is a homotopy $J:X\times\left[
0,1\right]  \rightarrow X$ such that

\begin{enumerate}
\item[i)] $J_{0}=id_{X}$,

\item[ii)] $\overline{J_{1}(X)}$ is compact, and

\item[iii)] $\left.  J\right\vert _{C\times\left[  0,1\right]  }%
=id_{C\times\left[  0,1\right]  }$.
\end{enumerate}

\begin{proof}
We need only prove the forward implication, as the converse is obvious. Begin
with a homotopy $H:X\times\left[  0,1\right]  \rightarrow X$ satisfying the
analogues of conditions i) and ii). Choose a compact polyhedral neighborhood
$D$ of $C$ in $X$. Then define $J$ on the union of $X\times\left\{  0\right\}
$ and $\left(  C\cup\left(  X-intD\right)  \right)  \times\left[  0,1\right]
$ as follows:%
\[
J\left(  x,t\right)  =\left\{
\begin{tabular}
[c]{rr}%
$x$ & if $t=0$\\
$x$ & if $x\in C$\\
$H\left(  x,t\right)  $ & if $x\in X-intD$%
\end{tabular}
\ \ \ \right.  \text{.}%
\]
Apply the Homotopy Extension Theorem \cite[\S IV.2]{Hu} to extend $J$ to all
of $X\times\left[  0,1\right]  $. Condition ii) follows from compactness of
$D$.
\end{proof}
\end{lemma}

If a space is finitely dominated, one often wishes to know whether it has
finite homotopy type. This issue was resolved by Wall in \cite{Wa} where, to
every finitely dominated $X$, there is defined an obstruction $\sigma\left(
X\right)  $ lying in the reduced projective class group $\widetilde{K}%
_{0}\left(  \mathbb{Z[}\pi_{1}(X)]\right)  $. This obstruction vanishes if and
only if $X$ has finite homotopy type.

A space having finite homotopy type may have infinitely many ends. One example
is the universal cover of a figure-eight. However, within the realm of open
manifolds, this does not happen. In fact, we have

\begin{lemma}
The number of ends of a finitely dominated open $n$-manifold $M^{n}$ is a
finite integer bounded above by $rank\left(  H_{n-1}\left(  M^{n}%
,\mathbb{Z}_{2}\right)  \right)  +1$.

\begin{proof}
[Sketch of Proof]If $M^{n}$ is dominated by a finite complex $L$, then each
homology group of $L$ surjects onto the corresponding homology group of
$M^{n}$. It follows that $rank\left(  H_{n-1}\left(  M^{n},\mathbb{Z}%
_{2}\right)  \right)  <\infty$. Next observe that, for an improved
neighborhood of infinity $N$, the collection of boundary components of $N$
forms a nearly independent collection of elements of $H_{n-1}\left(
M^{n},\mathbb{Z}_{2}\right)  $. (This is were we use the fact that $M^{n}$ is
an open manifold.) So if $M^{n}$ contained improved neighborhoods of infinity
with arbitrarily large numbers of components, $H_{n-1}\left(  M^{n}%
,\mathbb{Z}_{2}\right)  $ would be infinitely generated. See \cite[Prop.3.1]%
{GuTi} for details.
\end{proof}
\end{lemma}

\subsection{Fundamental group at infinity and Siebenmann's thesis}

For an inverse sequence
\[
G_{0}\overset{\lambda_{1}}{\longleftarrow}G_{1}\overset{\lambda_{2}%
}{\longleftarrow}G_{2}\overset{\lambda_{3}}{\longleftarrow}\cdots
\]
of groups and homomorphisms, a \emph{subsequence} of $\left\{  G_{i}%
,\lambda_{i}\right\}  $ is an inverse sequence of the form
\[
G_{i_{0}}\overset{\phi_{1}}{\longleftarrow}G_{i_{1}}\overset{\phi_{2}%
}{\longleftarrow}G_{i_{2}}\overset{\phi_{3}}{\longleftarrow}\cdots.
\]
where, for each $j$, the homomorphism $\phi_{j}$ is the obvious composition
$\lambda_{i_{j-1}+1}\circ\cdots\circ\lambda_{i_{j}}$ of homomorphisms from the
original sequence. We say that $\left\{  G_{i},\lambda_{i}\right\}  $ is
\emph{stable}, if it contains a subsequence $\left\{  G_{i_{j}},\phi
_{j}\right\}  $ that induces a sequence of isomorphisms%
\begin{equation}
im\left(  \phi_{1}\right)  \overset{\cong}{\longleftarrow}im\left(  \phi
_{2}\right)  \overset{\cong}{\longleftarrow}im\left(  \phi_{3}\right)
\overset{\cong}{\longleftarrow}\cdots. \tag{*}%
\end{equation}
If a sequence (*) exists where the bonding maps are simply injections, we say
that $\left\{  G_{i},\lambda_{i}\right\}  $ is \emph{pro-injective}; if one
exists where the bonding maps are surjections, we say that $\left\{
G_{i},\lambda_{i}\right\}  $ is \emph{pro-surjective} or (more commonly)
\emph{semistable.}

For a one-ended open manifold $M^{n}$ and a neat sequence $\left\{
N_{i}\right\}  _{i=1}^{\infty}$ of $0$-neighborhoods of infinity, choose
basepoints $p_{i}\in N_{i}$, and paths $\alpha_{i}\subset N_{i}$ connecting
$p_{i}$ to $p_{i+1}$. Then construct an inverse sequence of groups:
\begin{equation}
\pi_{1}\left(  N_{0},p_{0}\right)  \overset{\lambda_{1}}{\longleftarrow}%
\pi_{1}\left(  N_{1},p_{1}\right)  \overset{\lambda_{2}}{\longleftarrow}%
\pi_{1}\left(  N_{2},p_{2}\right)  \overset{\lambda_{3}}{\longleftarrow}%
\cdots. \tag{\dag}%
\end{equation}
by letting $\lambda_{i+1}:\pi_{1}\left(  N_{i+1},p_{i+1}\right)
\rightarrow\pi_{1}\left(  N_{i},p_{i}\right)  $ be the homomorphism induced by
inclusion followed by the change of basepoint isomorphism determined by
$\alpha_{i}$. The obvious singular ray obtained by piecing together the
$\alpha_{i}$'s is often referred to as the \emph{base ray }for this inverse
sequence. This inverse sequence (or more precisely the `pro-equivalence class'
of this sequence) is referred to as the \emph{fundamental group at infinity}
for $M^{n}$ and is denoted by $\pi_{1}\left(  \varepsilon\left(  M^{n}\right)
\right)  $.

\begin{remark}
For the purposes of this paper, we only need to consider the fundamental group
at infinity for \textbf{one-ended} manifolds. (Even though we often begin with
a multi-ended manifold.) In multi-ended situations, one may associate a
different inverse sequence to each end. For example, if $M^{n}$ is a $k$-ended
open manifold and $\left\{  N_{i}\right\}  _{i=1}^{\infty}$ as a well-indexed
neat sequence of\emph{ }$0$-neighborhoods of infinity. Then, for each
$j\in\left\{  1,2,\cdots,k\right\}  $ we can construct an inverse sequence
\[
\pi_{1}\left(  N_{0}^{j}\right)  \overset{\lambda_{1}^{j}}{\longleftarrow}%
\pi_{1}\left(  N_{1}^{j}\right)  \overset{\lambda_{2}^{j}}{\longleftarrow}%
\pi_{1}\left(  N_{2}^{j}\right)  \overset{\lambda_{3}^{j}}{\longleftarrow
}\cdots.
\]
which is called the \emph{fundamental group at the }$j^{th}$ \emph{end of
}$M^{n}$. (Here we have omitted reference to basepoints only to simplify notation.)
\end{remark}

For a more thorough discussion of inverse sequences and the fundamental group
system at infinity, see \cite{Gu}.\medskip

As indicated in the introduction, Theorem \ref{main theorem} will be obtained
as a consequence of the main result of \cite{Si}. For easy reference, we state
that result and provide some necessary definitions.

\begin{theorem}
[Siebenmann, 1965]\label{siebenmann}A one-ended open $n$-manifold $M^{n}$
($n\geq6$) is homeomorphic to the interior of a compact manifold with boundary iff:

\begin{enumerate}
\item $M^{n}$ is inward tame at infinity,

\item $\pi_{1}$ is stable at infinity, and

\item $\sigma_{\infty}\left(  M^{n}\right)  \in\widetilde{K}_{0}\left(
\mathbb{Z[}\pi_{1}(\varepsilon(M^{n}))]\right)  $ is trivial.\bigskip
\end{enumerate}
\end{theorem}

In the above, \emph{inward tame at infinity }(or simply `inward tame') means
that for any neighborhood $N$ of infinity, there exists a homotopy (sometimes
called a \emph{taming homotopy}) $H:N\times\left[  0,1\right]  \rightarrow V$
such that $H_{0}=id$ and $\overline{H_{1}(N)}$ is compact. Equivalently
$M^{n}$ is inward tame if all clean neighborhoods of infinity are finitely
dominated. If $N\supseteq N^{\prime}$ are clean neighborhoods of infinity,
then any taming homotopy for $N^{\prime}$ can be extended to a taming homotopy
for $N$. Thus, in order to prove inward tameness for $M^{n}$, it suffices to
show the existence of arbitrarily small finitely dominated clean neighborhoods
of infinity.

Given conditions 1) and 2) above, one may choose a $0$-neighborhood of
infinity $N$ with the `correct' fundamental group---as determined by 2). Then
$\sigma_{\infty}\left(  M^{n}\right)  $ is the Wall finiteness obstruction of
$N$. With some additional work, one sees that $\sigma_{\infty}\left(
M^{n}\right)  $ is trivial if and only if \emph{all} clean neighborhoods of
infinity in $M^{n}$ (or equivalently, arbitrarily small clean neighborhoods of
infinity) have finite homotopy type. For more details see \cite{Si} or
\cite{Gu}.

\begin{remark}
By giving a more general definition of $\sigma_{\infty}\left(  M^{n}\right)
$, it is possible to separate Conditions 2) and 3); this has been done in
\cite{Gu}. However, in the current context, it seems better to keep Theorem
\ref{siebenmann} in its traditional form.
\end{remark}

\subsection{Combinatorial group theory and the Generalized Seifert-VanKampen
Theorem}

The last bit of background information we wish to comment on is primarily
combinatorial group theory. Given groups $G_{0},G_{1}$ and $G_{2}$ and
homomorphisms $i_{1}:G_{0}\rightarrow G_{1}$ and $i_{2}:G_{0}\rightarrow
G_{2}$ we call $G$ the \emph{pushout} of $\left(  i_{1},i_{2}\right)  $ if
there exist homomorphisms $j_{1}$ and $j_{2}$ completing a commutative diagram%
\begin{equation}%
\begin{array}
[c]{ccc}%
G_{0} & \overset{i_{1}}{\rightarrow} & G_{1}\\
i_{2}\downarrow\quad &  & \quad\downarrow j_{2}\\
G_{2} & \overset{j_{1}}{\rightarrow} & G
\end{array}
\tag{$\diamondsuit$}%
\end{equation}
and satisfying the following `universal mapping property':%

\begin{gather*}
\text{\emph{If homomorphisms} }k_{1}:G_{1}\rightarrow H\text{ \emph{and}
}k_{2}:G_{2}\rightarrow H\text{ \emph{allow for a similar}}\\
\text{\emph{commutative diagram, then there exists a unique homomorphism}}\\
\varphi:G\rightarrow H\text{ \emph{such that} }j_{1}\circ\varphi=k_{1}\text{
\emph{and} }j_{2}\circ\varphi=k_{2}.
\end{gather*}
In this case, $G$ is uniquely determined up to isomorphism.

In the special case where the above homomorphisms $i_{1}$ and $i_{2}$ are
injective, the pushout is called a \emph{free product with amalgamation of
}$G_{1}$ and $G_{2}$ along $G_{0}$ and is denoted by $G_{1}\ast_{G_{0}}G_{2}$
With this terminology, we are implicitly viewing $G_{0}$ as a subgroup of both
$G_{1}$ and $G_{2}$. Then $G_{1}\ast_{G_{0}}G_{2}$ is the result of `gluing'
$G_{1}$ to $G_{2}$ along $G_{0}$. More precisely, if $\left\langle \left.
A_{1}\ \right\vert \ R_{1}\right\rangle $ and $\left\langle \left.
A_{1}\ \right\vert \ R_{1}\right\rangle $ are presentations for $G_{1}$ and
$G_{2}$ and $B$ is a generating set for $G_{0}$, then $\left\langle \left.
A_{1}\cup A_{2}\ \right\vert \ R_{1},R_{2},S\right\rangle $ is a presentation
for $G_{1}\ast_{G_{0}}G_{2}$ where%
\[
S=\left\{  \left.  i_{1}\left(  y\right)  ^{-1}i_{2}\left(  y\right)
\ \right\vert \ y\in B\right\}  .
\]
(The same procedure produces presentations for arbitrary pushouts.)\smallskip

In topology, the most common application of `pushout diagrams' is found in the
Seifert-VanKampen Theorem \cite[Ch.4]{Ma}, which may stated as follows: if a
space $X$ is expressed as a union $X=U\cup V$ of path connected open sets such
that $U\cap V$ is also path connected and $x\in U\cap V$, then $\pi_{1}\left(
X,x\right)  $ is the pushout of
\[%
\begin{array}
[c]{ccc}%
\pi_{1}\left(  U\cap V,x\right)  & \overset{\theta_{1}}{\rightarrow} & \pi
_{1}\left(  U,x\right) \\
\theta_{2}\downarrow &  & \\
\pi_{1}\left(  V,x\right)  &  &
\end{array}
\]
where $\theta_{1}$ and $\theta_{2}$ are induced by inclusion. In most cases,
$\theta_{1}$ and $\theta_{1}$ are not injective, so $\pi_{1}\left(
X,x\right)  $ is not necessarily a free product with amalgamation.

The above group theoretic constructions can be extended to \emph{generalized
graphs of groups} and the more restrictive \emph{graphs of groups. }For either
of these constructions, we begin with an oriented graph $\Delta$. (Here a
`graph' is simply a $1$-dimensional CW complex.) Then, to each vertex $v$ we
associate a `vertex group' $G_{v}$ and to each edge $e$ we associate an `edge
group' $G_{e}$. In addition, for each edge group $G_{e}$ we need `edge
homomorphisms' $\varphi_{e}^{+}$ and $\varphi_{e}^{-}$ mapping $G_{e}$ into
the vertex group at the `positive' and `negative' end of $e$, respectively.
(If $e$ is a loop in $\Delta$, then $\varphi_{e}^{+}$ and $\varphi_{e}^{-}$
can be different homomorphisms into the same vertex group. Let $\left(
\mathcal{G},\Delta\right)  $ represent this setup. If each edge homomorphism
is injective, call $\left(  \mathcal{G},\Delta\right)  $ a `graphs of groups';
otherwise it is just a `generalized graphs of groups'.

Our next task is to assign, to an arbitrary generalized graph of groups
$\left(  \mathcal{G},\Delta\right)  $, a single group that generalizes the
pushout construction for the simple case. This could be done with a universal
mapping property. Instead, we describe a specific construction of the group.
Let $V$ [resp., $E$] denote the collection of vertices [resp., edges] of
$\Delta$. Choose a maximal tree $\Upsilon$ in $\Delta$. Then the
\emph{fundamental group of }$\left(  \mathcal{G},\Delta\right)  $ \emph{based
at} $\Upsilon$ is the group%
\[
\pi_{1}\left(  \mathcal{G},\Delta;\Upsilon\right)  =((\ast_{v\in V}G_{v})\ast
F_{E})/N
\]
where $\ast_{v\in V}G_{v}$ is the free product of all vertex groups, $F_{E}$
is the free group generated by the set $E$ and $N$ is the smallest normal
subgroup of $(\ast_{v\in V}G_{v})\ast F_{E}$ generated by the set%
\[
\left\{  \left.  e^{-1}\cdot\varphi_{e}^{-}\left(  x\right)  \cdot
e\cdot\left(  \varphi_{e}^{+}\left(  x\right)  \right)  ^{-1}\ \right\vert
\ e\in E\text{ and }x\in G_{e}\right\}  \cup\left\{  \left.  e\ \right\vert
\ e\in E\right\}  .
\]

\begin{example}
Diagram ($\diamondsuit$) determines a generalized graph of groups where the
graph is simply an oriented interval; moreover, the fundamental group of that
generalized graph of groups is precisely the pushout of that diagram. When
$\iota_{1}$ and $i_{2}$ are injective we have a genuine graph of groups whose
fundamental group is a free product with amalgamation.

A similar special case---this time, a graph of groups with just one vertex and
one edge---leads to another well-known construction in combinatorial group
theory---the HNN extension.
\end{example}

See \cite{Co} for details on the above ideas.\smallskip

In group theory, it is preferable to study free products with amalgamation
over arbitrary pushouts. Similarly, genuine graphs of groups are preferable to
generalized graphs of groups. However, as noted above, arbitrary pushouts
occur naturally in topology via the classical Seifert-VanKampen Theorem.
Similarly, the following Generalized Seifert-VanKampen Theorem frequently
leads to a generalized graph of groups.

\begin{theorem}
[Generalized Seifert-VanKampen]Suppose a path connected space $X$ may be
expressed as a union of path connected open subsets $\left\{  U_{\alpha
}\right\}  _{\alpha\in A}$ so that no point of $X$ lies in more than two of
the $U_{\alpha}$'s. Let $\Delta$ be the graph having vertex set $A$, and one
edge between $U_{\alpha}$ and $U_{\alpha^{\prime}}$ for each path component
$V_{\alpha\alpha^{\prime}\beta}$ of $U_{\alpha}\cap U_{\alpha^{\prime}}$
($\alpha\neq\alpha^{\prime}\in A$). Place an arbitrarily chosen orientation on
each edge; then choose a base point from each $U_{\alpha}$ [resp., each
$V_{\alpha\alpha^{\prime}\beta}$] and associate to the corresponding vertex
[resp., edge] the fundamental group of $U_{\alpha}$ [resp., $V_{\alpha
\alpha^{\prime}\beta}$]. For each $V_{\alpha\alpha^{\prime}\beta}$, choose
paths in $U_{\alpha}$ and $U_{\alpha^{\prime}}$ respectively, connecting the
base point of $V_{\alpha\alpha^{\prime}\beta}$ to the base points of
$U_{\alpha}$ and $U_{\alpha^{\prime}}$. Let the two edge homomorphisms for
$V_{\alpha\alpha^{\prime}\beta}$ be those induced by inclusion followed by
change of basepoints along these paths. If $\left(  \mathcal{G},\Delta\right)
$ denotes this graph of groups, then $\pi_{1}(X,x)$ is isomorphic to $\pi
_{1}\left(  \mathcal{G},\Delta;\Upsilon\right)  $ where $\Upsilon$ is an
arbitrarily chosen maximal tree in $\Delta.$

\begin{proof}
See Chapter 2 of \cite{Ge}.
\end{proof}
\end{theorem}

\section{Proof of Theorem \ref{main theorem}\label{Section: Proof}}

In order to prove Theorem \ref{main theorem} we need only show that, if an
open manifold $M^{n}$ has finite homotopy type, then $M^{n}\times\mathbb{R}$
is a one-ended open manifold satisfying all three conditions of Theorem
\ref{siebenmann}. Proposition \ref{main-prop} begins that process; Part a)
asserts that $M^{n}\times\mathbb{R}$ is one-ended and open, while Part b)
ensures that Condition 1) holds. Strictly speaking, the `end obstruction'
$\sigma_{\infty}\left(  M^{n}\times\mathbb{R}\right)  $ found in Theorem
\ref{siebenmann} cannot be defined until it it known that Condition 2) holds.
Even so, it is possible to address Condition 3) before Condition 2) by proving
that \emph{all} clean neighborhoods of infinity in $M^{n}\times\mathbb{R}$
have finite homotopy type. This will be done (under the assumption that
$M^{n}$ has finite homotopy type) in Part c). Therefore, in the context of
Theorem \ref{main theorem}, once Condition 2) is verified, Condition 3)
follows immediately. To simplify the discussion, we refer to a manifold in
which all clean neighborhoods of infinity have finite homotopy type as
\emph{super-inward tame at infinity.}

\subsection{Conditions 1) and 3) of Theorem \ref{siebenmann}}

Before stating Proposition \ref{main-prop}, we introduce some terminology and
notation to used through the reest of this paper. Given a connected open
manifold $M^{n}$, a clean neighborhood of infinity $N\subseteq M^{n}$, and
$m>0$, the \emph{associated neighborhood of infinity} in $M^{n}\times
\mathbb{R}$ is the set
\[
W\left(  N,m\right)  =(N\times\mathbb{R)}\cup(M^{n}\times((-\infty
,-m]\cup\lbrack m,\infty)))
\]
It is easy to see that $W\left(  N,m\right)  $ is indeed a neighborhood of
infinity, that $W\left(  N,m\right)  $ is always connected, and that $W\left(
N,m\right)  $ is a $0$-neighborhood of infinity in $M^{n}\times\mathbb{R}$
whenever $N$ is a $0$-neighborhood of infinity $M^{n}$.

In addition, let
\begin{align*}
W^{+}\left(  N,m\right)   &  =(N\times\mathbb{R)}\cup(M^{n}\times\lbrack
m,\infty))\text{ and}\\
W^{-}\left(  N,m\right)   &  =(N\times\mathbb{R)}\cup(M^{n}\times(-\infty,-m])
\end{align*}
Then $W^{+}\left(  N,m\right)  $ deformation retracts onto $M^{n}%
\times\left\{  m\right\}  $ and $W^{-}\left(  N,m\right)  $ deformation
retracts onto $M^{n}\times\left\{  -m\right\}  $, so both are homotopy
equivalent to $M^{n}$. Moreover $W^{+}\left(  N,m\right)  $ $\cap W^{-}\left(
N,m\right)  =N\times\mathbb{R\ \simeq\ }N$.

\begin{proposition}
\label{main-prop}Let $M^{n}$ be a connected open $n$-manifold. Then

\begin{enumerate}
\item[a)] $M^{n}\times\mathbb{R}$ is a one-ended open $\left(  n+1\right)  $-manifold.

\item[b)] $M^{n}\times\mathbb{R}$ is inward tame at infinity iff $M^{n}$ is
finitely dominated.

\item[c)] $M^{n}\times\mathbb{R}$ is super-inward tame at infinity iff $M^{n}$
has finite homotopy type.
\end{enumerate}

\begin{proof}
As noted above, all neighborhoods of infinity of the type $W\left(
N,m\right)  $ are connected; moreover, they can be made arbitrarily small by
choosing $N$ to be small and $m$ large. Therefore, $M^{n}\times\mathbb{R}$ is one-ended.

The forward implications of Assertions a) and b) are immediate. In particular,
since $M^{n}\times\mathbb{R}$ itself is a clean neighborhood of infinity in
$M^{n}\times\mathbb{R}$, both implications can be deduced from the homotopy
equivalence $M^{n}\simeq M^{n}\times\mathbb{R}$. Thus, we turn our attention
to the two converses.

Given $W\left(  N,m\right)  $, let $C=M^{n}-int\left(  N\right)  $ (a compact
codimension $0$ submanifold of $M^{n}$). Then let $C^{\prime}$ denotes a
second `copy' of $C$ disjoint from $M^{n}$; and $K_{N}$ be the adjunction
space
\[
K_{N}=M^{n}\cup_{\varphi}C^{\prime}%
\]
obtained by attaching $C^{\prime}$ to $M^{n}$ along its boundary via the
`identity map' $\varphi:\partial C^{\prime}\rightarrow\partial C$.

It is easy to see that $K_{N}$ is homotopy equivalent to $W\left(  N,m\right)
$; indeed, $W\left(  N,m\right)  $ deformation retracts onto the subset%
\[
(M^{n}\times\left\{  m\right\}  )\cup\left(  \partial C\times\left[
-m,m\right]  \right)  \cup\left(  C\times\left\{  -m\right\}  \right)
\]
which is homeomorphic to $K_{N}$. Thus, to show that $M^{n}\times\mathbb{R}$
is inward tame at $\infty$, [resp., $M^{n}\times\mathbb{R}$ is super-inward
tame at $\infty$], it suffices to show that $K_{N}$ is finitely dominated
[resp., has finite homotopy type].

\begin{claim}
If $M^{n}$ is finitely dominated, then $K_{N}$ is finitely
dominated.\smallskip
\end{claim}

By Lemma \ref{improved-domination}, we may choose a homotopy $J:M^{n}%
\times\left[  0,1\right]  \rightarrow M^{n}$ such that $\left.  J\right\vert
_{(M^{n}\times\left\{  0\right\}  )\cup\left(  C\times\left[  0,1\right]
\right)  }$ is the identity, and $\overline{J_{1}\left(  M^{n}\right)  }$ is
compact. Extend $J$ to a homotopy $J^{\ast}:K_{N}\times\left[  0,1\right]
\rightarrow K_{V}$ by letting $J^{\ast}$ be the identity over $C^{\prime}$.
Then $J_{0}^{\ast}$ is the identity, and $J_{1}^{\ast}\left(  K_{V}\right)  $
has compact closure; so $K_{N}$ is finitely dominated.

\begin{claim}
If $M^{n}$ has finite homotopy type, then $K_{V}$ has finite homotopy
type.\smallskip
\end{claim}

Let $f:M^{n}\rightarrow L$ be a homotopy equivalence, where $L$ is a finite
complex. Then $K_{N}=M^{n}\cup_{\varphi}C^{\prime}$ is homotopy equivalent to
the adjunction space
\[
L\cup_{f\circ\varphi}C^{\prime}%
\]
where $f\circ\varphi$ maps $\partial C^{\prime}$ into $L$. This latter
adjunction space is homotopy equivalent to a finite complex. In fact, if we
begin with a triangulation of $M^{n}$ so that $C$ and $\partial C$ are
subcomplexes and choose $f$ to be a cellular map, then $f\circ\varphi$ is also
cellular and $L\cup_{f\circ\varphi}C^{\prime}$ is a finite complex.\medskip
\end{proof}
\end{proposition}

\subsection{Main Step: Stability of $\pi_{1}$ at infinity}

The following will show that $M^{n}\times\mathbb{R}$ satisfies Condition 2) of
Theorem \ref{siebenmann}, and thus complete our proof of Theorem
\ref{main theorem}.

\begin{proposition}
\label{pi1-stable}If a connected open manifold $M^{n}$ is finitely dominated,
then $M^{n}\times\mathbb{R}$ has stable fundamental group at infinity.
\end{proposition}

Let $M^{n}$ be a $k$-ended open $n$-manifold and $P$ and $Q$ be $0$%
-neighborhoods of infinity with $Q\subseteq int\left(  P\right)  $. Index the
components $P^{1},P^{2},\cdots,P^{k}$ of $P$ and $Q^{1},Q^{2},\cdots,Q^{k}$ of
$Q$ so that $Q^{j}\subseteq P^{j}$ for $j$ $=1,\cdots,k$. For each $j$ , let
$A^{j}=\overline{Q^{j}-P^{j}}$.

If $M^{n}$ is finitely dominated, choose $P$ sufficiently small that there is
a homotopy $H:M^{n}\times\left[  0,1\right]  \rightarrow M^{n}$ pulling
$M^{n}$ into $M^{n}-P$. In addition (by Lemma \ref{improved-domination})
arrange that $H$ is fixed over some non-empty open set $U$. To simplify
notation, we focus on a single end; in particular, let $j\in\left\{
1,\cdots,k\right\}  $ be fixed. Choose basepoints $p_{\ast}\in U$,
$p\in\partial P^{j}$ and $q\in\partial Q^{j}$; then choose a proper embedding
$r:[0,\infty)\rightarrow M^{n}$ such that $r\left(  0\right)  =p_{\ast}$,
$r\left(  1\right)  =p$, $r\left(  2\right)  =q$, and so that the image ray
$R=r\left(  [0,\infty\right)  )$ intersects each of $\partial P^{j}$ and
$\partial Q^{j}$ transversely once and only once at the points $p$ and $q$,
respectively. Let $\alpha=R\cap A^{j}$ denote the corresponding arc in $A^{j}$
between $p$ and $q$.

Let $t:B^{n-1}\times\lbrack-1,\infty)\rightarrow M^{n}$ be a homeomorphism
onto a regular neighborhood $T$ of $R$ so that $\left.  t\right\vert
_{\left\{  \overline{0}\right\}  \times\lbrack0,\infty)}=r$, and so that
$T\cap A^{j}$ is a relative regular neighborhood of $\alpha$ in $A^{j}$
intersecting $\partial P^{j}$ and $\partial Q^{j}$ in $\left(  n-1\right)
$-disks $D$ and $D^{\prime}$, with $D=t\left(  B^{n-1}\times\left\{
1\right\}  \right)  $ and $D^{\prime}=t\left(  B^{n-1}\times\left\{
2\right\}  \right)  $. Then choose an $\left(  n-1\right)  $-ball
$B_{0}\subseteq intB^{n-1}$, centered at $\overline{0}$; and let
$T_{0}=t\left(  B_{0}\times\lbrack-1,\infty)\right)  $ be the corresponding
smaller regular neighborhood of $R$, with corresponding subdisks $D_{0}$ and
$D_{0}^{\prime}$ contained in $intD$ and $intD^{\prime}$, respectively. We now
utilize the `homotopy refinement procedure' developed on pages 267-268 of
\cite{GuTi} to replace $H$ with a new homotopy $K:M^{n}\times\left[
0,1\right]  \rightarrow M^{n}$ which, in addition to pulling $M^{n}$ into
$M^{n}-P$, has the properties:

\begin{itemize}
\item[i)] $K$ is `canonical' over $T_{0}$, and

\item[ii)] tracks of points lying outside $T_{0}$ do not pass through the
interior of $T_{0}$.
\end{itemize}

The first of these properties arranges that all tracks of points in $R$
proceed monotonically in $R$ to $p$; and that $\left.  K\right\vert
_{D_{0}^{\prime}\times\left[  0,1\right]  }$ takes \newline$D_{0}^{\prime
}\times\left[  0,\frac{1}{2}\right]  $ homeomorphically onto $t\left(
B_{0}\times\left[  0,2\right]  \right)  $, with $D_{0}^{\prime}\times\left[
\frac{1}{2},1\right]  $ being flattened onto $t\left(  B_{0}\times\left\{
0\right\}  \right)  $. We may also arrange that $K\left(  D_{0}^{\prime}%
\times\left\{  \frac{1}{4}\right\}  \right)  =D_{0}$. See \cite{GuTi} for details.

\begin{proposition}
\label{pushing-loops}Assume the above setup, with $M^{n}$ finitely dominated,
$j\in\left\{  1,\cdots,k\right\}  $ fixed, and all previous notation
unchanged. Then every loop $\tau$ in $P^{j}$ based at $p$ is homotopic (rel
$p$) in $M^{n}$ to a loop of the form $\alpha\ast\tau^{\prime}\ast\alpha^{-1}%
$, where $\tau^{\prime}$is a loop in $Q^{j}$ based at $q$.

\begin{proof}
Every loop in $P^{j}$ based at $p$ is homotopic (rel $p$) to a product
$\tau_{1}\ast\tau_{2}\ast\cdots\tau_{u}$ where: for each $v$, either $\tau
_{v}$ lies entirely in $A^{j}$ or $\tau_{v}$ is (already) a loop of the form
$\alpha\ast\tau_{v}^{\prime}\ast\alpha^{-1}$ where $\tau_{v}^{\prime}$ is a
loop in $Q^{j}$ based at $q$. So, without loss of generality, we assume, that
$\tau$ lies entirely in $A^{j}$.

Consider the map $L=\left.  K\right\vert _{\partial Q^{j}\times\left[
0,1\right]  }:\partial Q^{j}\times\left[  0,1\right]  \rightarrow M^{n}$.
Choose triangulations $\Delta_{1}$ and $\Delta_{2}$ of the domain and range,
respectively. Without changing its definition on $(\partial Q^{j}%
\times\left\{  0\right\}  )\cup(D_{0}^{\prime}\times\left[  0,\frac{1}%
{2}\right]  )$, adjust $L$ (up to a small homotopy) to a non-degenerate
simplicial map. Then adjust $\tau$ (rel $p$) to an embedded circle in general
position with respect to $\Delta_{2}$, lying entirely in $intA^{j}$, except at
its basepoint $p$, which lies in $\partial A^{j}$. Then $L^{-1}\left(
\tau\right)  $ is a closed $1$-manifold lying in $\partial Q^{j}\times(0,1)$.
Let $\sigma$ be the component of $L^{-1}\left(  \tau\right)  $ containing the
point $\left(  q,\frac{1}{4}\right)  $. Since $L$ takes a neighborhood of
$\left(  q,\frac{1}{4}\right)  $ homeomorphically onto a neighborhood of $p$,
and since no other points of $\sigma$ are taken near $p$ (use property ii)
above), then $L$ takes $\sigma$ onto $\tau$ in a degree $1$ fashion. Now the
natural deformation retraction of $\partial Q\times\left[  0,1\right]  $ onto
$\partial Q\times\left\{  0\right\}  $ pushes $\sigma$ into $\partial
Q\times\left\{  0\right\}  $, while sliding $\left(  q,\frac{1}{4}\right)  $
along the arc $\left\{  q\right\}  \times\left[  0,\frac{1}{4}\right]  $.
Composing this push with $L$ provides a homotopy of $\tau$ to a loop
$\tau^{\prime}$ in $\partial Q$ whereby, the basepoint $p$ is slid along
$\alpha$ to $q$. This provides the desired (basepoint preserving) homotopy
from $\tau$ to $\alpha\ast\tau^{\prime}\ast\alpha^{-1}$.
\end{proof}
\end{proposition}

\begin{corollary}
\label{isomorphic images}Assume the full setup for Proposition
\ref{pushing-loops} and let
\begin{align*}
\Gamma_{P^{j}}  &  =im\left(  \pi_{1}\left(  P^{j},p\right)  \rightarrow
\pi_{1}\left(  M^{n},p\right)  \right)  \text{, and}\\
\Gamma_{Q^{j}}  &  =im\left(  \pi_{1}\left(  Q^{j},q\right)  \rightarrow
\pi_{1}\left(  M^{n},q\right)  \right)  .
\end{align*}
Then the change of basepoint isomorphism $\widehat{\alpha}:\pi_{1}\left(
M^{n},q\right)  \rightarrow\pi_{1}\left(  M^{n},p\right)  $ takes
$\Gamma_{Q^{j}}$ isomorphically onto $\Gamma_{P^{j}}$.

\begin{proof}
Since $\alpha\cup Q^{j}\subseteq P^{j}$, it is clear that $\widehat{\alpha}$
takes $\Gamma_{Q^{j}}$ into $\Gamma_{P^{j}}$. Injectivity is immediate, and
Proposition \ref{pushing-loops} assures surjectivity.\medskip
\end{proof}
\end{corollary}

We now turn our attention back to the manifold $M^{n}\times\mathbb{R}$. In
order to understand the fundamental group system at infinity, it will suffice
to understand `special' neighborhoods of infinity of the sort $W\left(
N,m\right)  $ (along with corresponding bonding maps). To simplify the
exposition, we first consider the special case where $M^{n}$ itself is
one-ended. Afterwards we upgrade the proof so that it includes the general case.

\begin{proposition}
\label{technical prop-one-ended}Suppose $M^{n}$ is a one-ended open
$n$-manifold and $P$ and $Q$ are $0$-neighborhoods of infinity in $M^{n}$ with
$Q\subseteq intP$. Choose $p\in\partial P$, $q\in\partial Q$ and $\alpha$ a
path in $\overline{P-Q}$ connecting $p$ to $q$. For $0<m<m^{\prime}<\infty$,
let $W\left(  P,m\right)  \supseteq W\left(  Q,m^{\prime}\right)  $ be
corresponding neighborhoods of infinity in $M^{n}\times\mathbb{R}$ and
$\lambda:\pi_{1}\left(  W\left(  Q,m^{\prime}\right)  ,(q,0)\right)
\rightarrow\pi_{1}\left(  W\left(  P,m\right)  ,(p,0)\right)  $ the
homomorphism induced by inclusion followed by a change of basepoints along
$\alpha\times0$. Then

\begin{enumerate}
\item $\pi_{1}\left(  W\left(  P,m\right)  ,(p,0)\right)  \cong\pi_{1}\left(
M^{n},p\right)  \ast_{\Gamma_{P}}\pi_{1}\left(  M^{n},p\right)  $, where
\[
\Gamma_{P}=im\left(  \pi_{1}\left(  P,p\right)  \rightarrow\pi_{1}\left(
M^{n},p\right)  \right)  ,
\]

\item $\pi_{1}\left(  W\left(  Q,m^{\prime}\right)  ,(q,0)\right)  \cong
\pi_{1}\left(  M^{n},q\right)  \ast_{\Gamma_{Q}}\pi_{1}\left(  M^{n},q\right)
$, where
\[
\Gamma_{Q}=im\left(  \pi_{1}\left(  Q,q\right)  \rightarrow\pi_{1}\left(
M^{n},q\right)  \right)  ,
\]

\item the homomorphism $\lambda$ is surjective, and

\item if there exists a homotopy pulling $M^{n}$ into $M^{n}-P$, then
$\lambda$ is an isomorphism.
\end{enumerate}

\begin{proof}
Using our earlier notation, the Seifert-VanKampen Theorem establishes
\newline$\pi_{1}\left(  W\left(  P,m\right)  ,(p,0)\right)  $ as the pushout
of the diagram\smallskip%
\[%
\begin{array}
[c]{ccc}%
\pi_{1}\left(  P\times\mathbb{R},(p,0)\right)   & \rightarrow\smallskip &
\pi_{1}\left(  W^{+}\left(  P,m\right)  ,(p,0)\right)  \\
\downarrow\smallskip &  & \\
\pi_{1}\left(  W^{-}\left(  P,m\right)  ,(p,0)\right)   &  &
\end{array}
\smallskip
\]
where both homomorphisms are induced by inclusion. Homotopy equivalences
\begin{align*}
\left(  P\times\mathbb{R},(p,0)\right)   &  \simeq\left(  P,p\right)  \text{,
and}\\
\left(  W^{+}\left(  P,m\right)  ,(p,0)\right)   &  \simeq\left(
M^{n},p\right)  \simeq\left(  W^{-}\left(  P,m\right)  ,(p,0)\right)
\end{align*}
allow us to replace the above with a simpler diagram%
\[%
\begin{array}
[c]{ccc}%
\pi_{1}\left(  P,p\right)  \smallskip & \overset{i_{\ast}}{\rightarrow} &
\pi_{1}\left(  M^{n},p\right)  \\
i_{\ast}\downarrow\smallskip &  & \\
\pi_{1}\left(  M^{n},p\right)   &  &
\end{array}
.
\]
This diagram does not define a free product with amalgamation since $i_{\ast}$
needn't be injective, however, the pushout is identical to that of
\begin{equation}%
\begin{array}
[c]{ccc}%
\Gamma_{P}\smallskip & \rightarrow & \pi_{1}\left(  M^{n},p\right)  \\
\downarrow\smallskip &  & \\
\pi_{1}\left(  M^{n},p\right)   &  &
\end{array}
\tag{\#}%
\end{equation}
(both homomorphisms are inclusions) which determines the free product with
amalgamation promised in 1).

Of course, assertion 2) is identical to 1). Then, from (\#) and the analogous
diagram for $\pi_{1}\left(  W\left(  Q,m^{\prime}\right)  ,(q,0)\right)  $, we
see that $\lambda$ is induced by the isomorphism $\pi_{1}\left(
M^{n},q\right)  \ast\pi_{1}\left(  M^{n},q\right)  \rightarrow\pi_{1}\left(
M^{n},p\right)  \ast\pi_{1}\left(  M^{n},p\right)  $ by quotienting out (in
the domain and range) by relations induced by $\Gamma_{Q}$ and $\Gamma_{P}$,
respectively---according to the definition of free product with amalgamation.
Thus, $\lambda$ is necessarily surjective. Moreover, if there is a homotopy
pulling $M^{n}$ into $P$, Corollary \ref{isomorphic images} ensures that
$\lambda$ is an isomorphism.
\end{proof}
\end{proposition}

\begin{corollary}
Let $M^{n}$ be a connected, one-ended open $n$-manifold. Then $M^{n}%
\times\mathbb{R}$ is a one-ended open $\left(  n+1\right)  $-manifold with
semistable fundamental group at infinity. If $M^{n}$ is finitely dominated,
then $M^{n}\times\mathbb{R}$ has stable fundamental group at infinity which is
pro-isomorphic to the system%
\[
\pi_{1}\left(  M^{n}\right)  \ast_{\Gamma}\pi_{1}\left(  M^{n}\right)
\overset{id}{\longleftarrow}\pi_{1}\left(  M^{n}\right)  \ast_{\Gamma}\pi
_{1}\left(  M^{n}\right)  \overset{id}{\longleftarrow}\pi_{1}\left(
M^{n}\right)  \ast_{\Gamma}\pi_{1}\left(  M^{n}\right)  \overset
{id}{\longleftarrow}\cdots
\]
where $\Gamma$ is the image (translated by an appropriate change of basepoint
isomorphism) of the fundamental group of any sufficiently small $0$%
-neighborhood of infinity in $\pi_{1}\left(  M^{n}\right)  $.

\begin{proof}
This corollary is almost immediate. One simply chooses a neat sequence
$\left\{  N_{i}\right\}  _{i=1}^{\infty}$ of $0$-neighborhoods of infinity in
$M^{n}$, then applies the previous proposition (repeatedly) to the sequence
$\left\{  W\left(  N_{i},i\right)  \right\}  _{i=1}^{\infty}$. If $M^{n}$ is
finitely dominated, $N_{1}$ should chosen sufficiently small that $M^{n}$ can
be pulled into $M^{n}-N_{1}$.
\end{proof}
\end{corollary}

In the introduction, we noted that, without the hypothesis of finite
domination on $M^{n}$, the fundamental group at infinity in $M^{n}%
\times\mathbb{R}$ needn't be stable. This is now easy to exhibit.

\begin{example}
[An $M^{n}\times\mathbb{R}$ with nonstable $\pi_{1}$ at infinity]Let
$T_{1}=B^{n-1}\times S^{1}\subseteq S^{n-1}\times S^{1}$ where $B^{n-1}%
\subseteq S^{n-1}$ is a tamely embedded $\left(  n-1\right)  $-ball. Then let
$T_{2}\subseteq int(T_{1})$ be another (thinner) copy of $B^{n-1}\times S^{1}$
that winds around twice in the $S^{1}$ direction. Inside $T_{2}$ choose a
third (even thinner) copy of $B^{n-1}\times S^{1}$ that winds through $T_{2}$
twice in the $S^{1}$ direction---and thus, four times through $T_{1}$ in the
original $S^{1}$ direction. Continue this infinitely to get a nested sequence
$T_{1}\supseteq T_{2}\supseteq T_{3}\supseteq\cdots$ so that $T_{\infty
}\subseteq S^{n-1}\times S^{1}$ is the \emph{dyadic solenoid}. Then
$M^{n}=(S^{n-1}\times S^{1})-T_{\infty}$ is a one-ended open $n$-manifold and
each $N_{i}=T_{i}-T_{\infty}$ is a $0$-neighborhood of infinity. Provided
$n>3$, it is easy to see that $N_{i}\hookrightarrow T_{i}$ induces a $\pi_{1}%
$-isomorphism. Therefore, the inverse sequence
\[
\pi_{1}\left(  N_{0},p_{0}\right)  \overset{\lambda_{1}}{\longleftarrow}%
\pi_{1}\left(  N_{1},p_{1}\right)  \overset{\lambda_{2}}{\longleftarrow}%
\pi_{1}\left(  N_{2},p_{2}\right)  \overset{\lambda_{3}}{\longleftarrow}%
\cdots.
\]
is isomorphic to the sequence
\[
\mathbb{Z}\overset{\times2}{\longleftarrow}\mathbb{Z}\overset{\times
2}{\longleftarrow}\mathbb{Z}\overset{\times2}{\longleftarrow}\mathbb{Z}%
\overset{\times2}{\longleftarrow}\mathbb{\cdots}\text{.}%
\]
A more descriptive form of the above inverse sequence is:%
\[
\mathbb{Z}\hookleftarrow2\mathbb{Z}\hookleftarrow4\mathbb{Z}\hookleftarrow
8\mathbb{Z\hookleftarrow\cdots}\text{.}%
\]
It follows that, for the corresponding sequence $\left\{  W\left(
N_{i},i\right)  \right\}  _{i=1}^{\infty}$ of neighborhoods of infinity in
$M^{n}\times\mathbb{R}$ we obtain a representation of the fundamental group at
infinity for $M^{n}\times\mathbb{R}$ isomorphic to the sequence%
\[
\mathbb{Z\ast}_{\mathbb{Z}}\mathbb{Z\leftarrow Z\ast}_{2\mathbb{Z}%
}\mathbb{Z\leftarrow Z\ast}_{4\mathbb{Z}}\mathbb{Z\leftarrow Z\ast
}_{8\mathbb{Z}}\mathbb{Z\leftarrow\cdots},
\]
where each bonding map is induced by the identity $\mathbb{Z\ast Z\rightarrow
Z\ast Z}$. Thus, each bond is surjective but not injective. It is easy to see
that such a system cannot be stable.
\end{example}

\begin{remark}
It is interesting to see that `crossing with $\mathbb{R}$' takes examples with
nonstable but pro-injective fundamental groups at infinity and produces
examples with nonstable but pro-surjective (semistable) fundamental groups at
infinity.\medskip
\end{remark}

We are now prepared to address the general situation where $M^{n}$ is
$k$-ended ($1\leq k<\infty$). For $k>1$, calculation of $\pi_{1}\left(
W\left(  N,m\right)  \right)  $ is complicated by the fact that $W^{+}\left(
N,m\right)  \cap W^{-}\left(  N,m\right)  =N\times\mathbb{R}$ is not
connected. In this situation, $\pi_{1}\left(  W\left(  N,m\right)  \right)  $
is most effectively described using a (generalized or actual) graph of groups.
In particular, let $\Theta_{k}$ denote the oriented graph consisting of two
vertices $v^{+}$ and $v^{-}$ and $k$ oriented edges $e^{1},e^{2},\cdots,e^{k}$
each running from $v^{-}$ to $v^{+}$. If $N^{1},N^{2},\cdots,N^{k}$ are the
components of $N$ with basepoints $p^{1}$ $,p^{2}$ $,\cdots,p^{k}$
respectively, we associate the following groups and homomorphisms to
$\Theta_{k}$:

\begin{itemize}
\item $G\left(  v^{+}\right)  =\pi_{1}\left(  W^{+}\left(  N,m\right)
,\left(  p^{1},0\right)  \right)  $ and $G\left(  v^{-}\right)  =\pi
_{1}\left(  W^{-}\left(  N,m\right)  ,\left(  p^{1},0\right)  \right)  ,$

\item for each $j\in\left\{  1,2,\cdots,k\right\}  $, $G\left(  e^{j}\right)
=\pi_{1}\left(  N^{j}\times\mathbb{R},\left(  p^{j},0\right)  \right)  ,$

\item For each $j\in\left\{  1,2,\cdots,k\right\}  $ the homomorphism
$\varphi_{j}^{+}:G\left(  e^{j}\right)  \rightarrow G\left(  v^{+}\right)  $
is the composition%
\[
\pi_{1}\left(  N^{j}\times\mathbb{R},\left(  p^{j},0\right)  \right)
\overset{i_{\ast}}{\longrightarrow}\pi_{1}\left(  W^{+}\left(  N,m\right)
,\left(  p^{j},0\right)  \right)  \overset{\widehat{\beta^{j}}}%
{\longrightarrow}\pi_{1}\left(  W^{+}\left(  N,m\right)  ,\left(
p^{1},0\right)  \right)
\]
where $\beta^{j}$ is an appropriately chosen path in $W^{+}\left(  N,m\right)
$ from $\left(  p^{j},0\right)  $ to $\left(  p^{1},0\right)  $.

\item The homomorphisms $\varphi_{j}^{-}:G\left(  e^{j}\right)  \rightarrow
G\left(  v^{-}\right)  $ are defined similarly to the above, but with $\beta$
chosen to lie in $W^{-}\left(  N,m\right)  $.
\end{itemize}

\noindent Since $\varphi_{j}^{+}$ and $\varphi_{j}^{-}$ needn't be injective,
the above setup is just a generalized\emph{ }graph of groups. Let it be
denoted by $(\mathcal{G}\left(  N\right)  ,\Theta_{k})$

Note that the edge $e^{1}$, by itself, is a maximal tree in $\Theta_{k}$. By
the Generalized Seifert-VanKampen Theorem $\pi_{1}\left(  N,p^{1}\right)
\cong\pi_{1}\left(  \mathcal{G}\left(  N\right)  ,\Theta_{k};e^{1}\right)  $.

We may obtain a similar---but simpler---graph of groups as follows. Again we
start with the graph $\Theta_{k}$. Motivated by the homotopy equivalences
$W^{+}\left(  N,m\right)  \simeq M^{n}\simeq W^{-}\left(  N,m\right)  $,
define both $G^{\prime}\left(  v^{+}\right)  $ and $G^{\prime}\left(
v^{+}\right)  $ to be $\pi_{1}\left(  M^{n},p^{1}\right)  $. Then, in order to
obtain injective edge homomorphisms, for each $j\in\left\{  1,2,\cdots
,k\right\}  $, let
\[
G^{\prime}\left(  e^{j}\right)  =\Gamma_{N^{j}}=im\left(  \pi_{1}\left(
N^{j},p^{j}\right)  \overset{i_{\ast}}{\longrightarrow}\pi_{1}\left(
M^{n},p^{j}\right)  \overset{\widehat{\beta^{j\prime}}}{\longrightarrow}%
\pi_{1}\left(  M^{n},p^{1}\right)  \right)  .
\]
and let all edge homomorphisms be inclusions. Here $\beta^{j\prime}$ is a path
in $M^{n}$ `parallel' to the path $\beta^{j}$ used above. This new graph of
groups will be denoted $\left(  \mathcal{G}^{\prime}\left(  N\right)
,\Theta_{k}\right)  $. It is easy to see a canonical isomorphism between
$\pi_{1}\left(  \mathcal{G}\left(  N\right)  ,\Theta_{k};e^{1}\right)  $ and
$\pi_{1}\left(  \mathcal{G}^{\prime}\left(  N\right)  ,\Theta_{k}%
;e^{1}\right)  $.

We are now ready to state a more general version of Proposition
\ref{technical prop-one-ended}, suitable for multi-ended $M^{n}$.

\begin{proposition}
\label{technical prop-k-ended}Suppose $M^{n}$ is a $k$-ended open $n$-manifold
and $P$ and $Q$ are $0$-neighborhoods of infinity in $M^{n}$ with components
of $P^{1},P^{2},\cdots,P^{k}$ and $Q^{1},Q^{2},\cdots,Q^{k}$, such that
$Q^{j}\subseteq int(P^{j})$ for each $j$. Choose basepoints $p^{j}$
$\in\partial P^{j}$ and $q_{j}\in\partial Q^{j}$, paths $\alpha^{j}$ in
$\overline{P^{j}-Q^{j}}$ connecting $p^{j}$ to $q^{j}$ and paths $\beta^{j}$
in $M^{n}$ connecting $q_{j}$ to $q_{1}$ for each $j\in\left\{  1,2,\cdots
,k\right\}  $. For $0<m<m^{\prime}<\infty$, let $W\left(  P,m\right)
\supseteq W\left(  Q,m^{\prime}\right)  $ be corresponding neighborhoods of
infinity in $M^{n}\times\mathbb{R}$, and $\lambda:\pi_{1}\left(  W\left(
Q,m^{\prime}\right)  ,(q_{1},0)\right)  \rightarrow\pi_{1}\left(  W\left(
P,m\right)  ,(p^{1},0)\right)  $ the homomorphism induced by inclusion
followed by a change of basepoints along $\alpha^{1}\times0$. Then

\begin{enumerate}
\item $\pi_{1}\left(  W\left(  P,m\right)  ,(p^{1},0)\right)  \cong\pi
_{1}\left(  \mathcal{G}^{\prime}\left(  P\right)  ,\Theta_{k};e^{1}\right)  $,
where $\left(  \mathcal{G}^{\prime}\left(  N\right)  ,\Theta_{k}\right)  $ is
the graph of groups described below with%
\[
\Gamma_{P^{j}}=im\left(  \pi_{1}\left(  P^{j},p^{j}\right)  \overset{i_{\ast}%
}{\longrightarrow}\pi_{1}\left(  M^{n},p^{j}\right)  \overset{\widehat
{\gamma^{j}}}{\longrightarrow}\pi_{1}\left(  M^{n},p^{1}\right)  \right)  ,
\]
and $\gamma^{j}=\alpha^{j}\ast\beta^{j}\ast(\alpha^{1})^{-1}$

\item $\pi_{1}\left(  W\left(  Q,m^{\prime}\right)  ,(q_{1},0)\right)
\cong\pi_{1}\left(  \mathcal{G}^{\prime}\left(  Q\right)  ,\Theta_{k}%
;e^{1}\right)  $ ($\left(  \mathcal{G}^{\prime}\left(  Q\right)  ,\Theta
_{k}\right)  $ analogous to figure below) where
\[
\Gamma_{Q^{j}}=im\left(  \pi_{1}\left(  Q^{j},q_{j}\right)  \overset{i_{\ast}%
}{\longrightarrow}\pi_{1}\left(  M^{n},q_{j}\right)  \overset{\widehat
{\beta^{j}}}{\longrightarrow}\pi_{1}\left(  M^{n},q_{1}\right)  \right)  ,
\]

\item the homomorphism $\lambda$ is surjective, and

\item if there exists a homotopy pulling $M^{n}$ into $P$, then $\lambda$ is
an isomorphism.
\end{enumerate}

\hspace*{1in}%
{\parbox[b]{3.8052in}{\begin{center}
\includegraphics[
trim=0.000000in -0.117993in 0.000000in 0.000000in,
height=1.4062in,
width=3.8052in
]%
{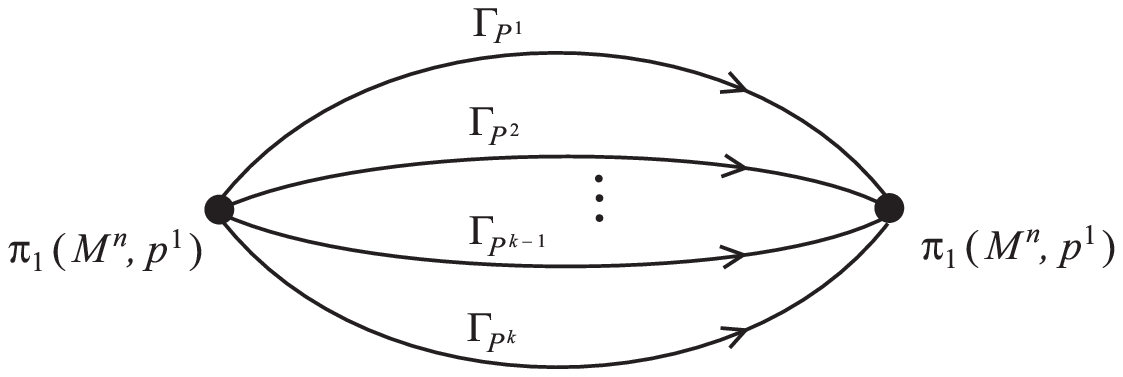}%
\\
$\left(  \mathcal{G}^{\prime}\left(  P\right)  ,\Theta_{k}\right)  $%
\end{center}}}%

\begin{proof}
As noted in the comments preceding this Proposition, 1) and 2) are essentially
just applications of the Generalized Seifert-VanKampen Theorem. Assertion 3)
is valid for nearly the same reason as Assertion 3) of Proposition
\ref{technical prop-one-ended}; in this case, $\lambda:\pi_{1}\left(
\mathcal{G}^{\prime}\left(  Q\right)  ,\Theta_{k};e^{1}\right)  \rightarrow
\pi_{1}\left(  \mathcal{G}^{\prime}\left(  P\right)  ,\Theta_{k};e^{1}\right)
$ is induced by the natural isomorphism
\[
\pi_{1}\left(  M^{n},q\right)  \ast\pi_{1}\left(  M^{n},q\right)  \ast
F_{E}\rightarrow\pi_{1}\left(  M^{n},p\right)  \ast\pi_{1}\left(
M^{n},p\right)  \ast F_{E}%
\]
where $F_{E}$ is the free group on generators $E=\left\{  e^{1},e^{2}%
,\cdots,e^{k}\right\}  $, by taking appropriate quotients in the domain and
the range (as prescribed by the definition of the fundamental group of a graph
of groups). If there exists a homotopy pulling $M^{n}$ into $P$, Corollary
\ref{pushing-loops} makes it clear that this homomorphism is an isomorphism.
\end{proof}
\end{proposition}

\begin{corollary}
\label{Corollary-k-ended}Let $M^{n}$ be a connected, $k$-ended open
$n$-manifold. Then $M^{n}\times\mathbb{R}$ is a one-ended open $\left(
n+1\right)  $-manifold with semi-stable fundamental group at infinity. If
$M^{n}$ is finitely dominated, then $M^{n}\times\mathbb{R}$ has stable
fundamental group at infinity which is pro-isomorphic to the system%
\[
\pi_{1}\left(  \mathcal{G}^{\prime},\Theta_{k};e^{1}\right)  \overset
{id}{\longleftarrow}\pi_{1}\left(  \mathcal{G}^{\prime},\Theta_{k}%
;e^{1}\right)  \overset{id}{\longleftarrow}\pi_{1}\left(  \mathcal{G}^{\prime
},\Theta_{k};e^{1}\right)  \overset{id}{\longleftarrow}\cdots
\]
where $\left(  \mathcal{G}^{\prime},\Theta_{k}\right)  $ is the graph of
groups pictured below. Here each $\Gamma_{j}$ is the image---translated by an
appropriate change of basepoint isomorphism---of the fundamental group of the
$j^{th}$ component of any sufficiently small $0$-neighborhood of infinity in
$\pi_{1}\left(  M^{n}\right)  $; and the edge homomorphisms are all inclusions.

\hspace*{1in}%
{\parbox[b]{3.4592in}{\begin{center}
\includegraphics[
trim=0.000000in -0.118065in 0.000000in 0.000000in,
height=1.4183in,
width=3.4592in
]%
{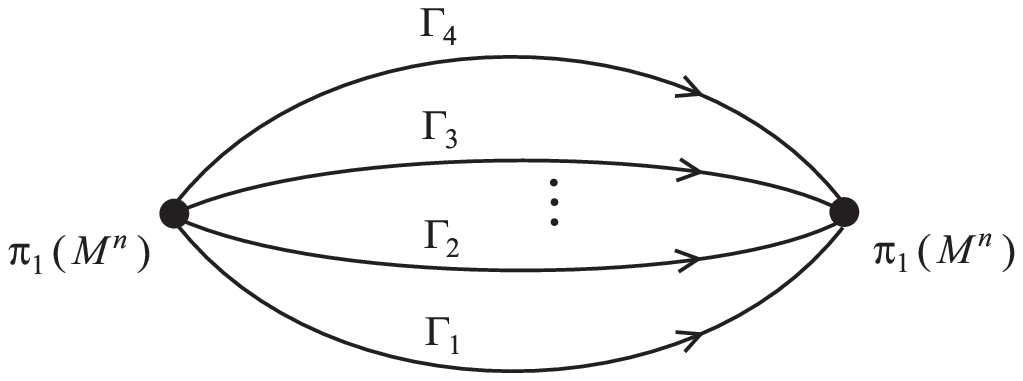}%
\\
$\left(  \mathcal{G}^{\prime},\Theta_{k}\right)  $%
\end{center}}}%

\begin{proof}
Choose a well-indexed neat sequence\emph{ }$\left\{  N_{i}\right\}
_{i=1}^{\infty}$ of $0$-neighborhoods of infinity in $M^{n}$. Then apply the
above Proposition repeatedly to the sequence $\left\{  W\left(  N_{i}%
,i\right)  \right\}  _{i=1}^{\infty}$ of associated neighborhoods of infinity
in $M^{n}\times\mathbb{R}$. If $M^{n}$ is finitely dominated, $N_{1}$ should
chosen sufficiently small that $M^{n}$ can be pulled into $M^{n}-N_{1}$.
\end{proof}
\end{corollary}

\section{Closing comments}

As indicated in the `easy part' of Theorem \ref{main theorem}, if $M^{n}$ is a
finitely dominated open manifold that is not homotopy equivalent to a finite
complex, then $M^{n}\times\mathbb{R}$ is not compactifiable by the addition of
a manifold boundary. My comparing Propositions \ref{main-prop} and
\ref{pi1-stable} with Theorem \ref{siebenmann}, it must be the case that
$\sigma_{\infty}\left(  M^{n}\times\mathbb{R}\right)  $ is non-trivial. Since
$M^{n}$ does not have finite homotopy type, its Wall obstruction
$\sigma\left(  M^{n}\right)  $ is a nontrivial element of $\widetilde{K}%
_{0}\left(  \mathbb{Z[}\pi_{1}(M^{n})]\right)  $. As one might expect, there
is a relationship between $\sigma_{\infty}\left(  M^{n}\times\mathbb{R}%
\right)  $ and $\sigma\left(  M^{n}\right)  $.

By Proposition \ref{technical prop-k-ended} and Corollary
\ref{Corollary-k-ended}, $\sigma_{\infty}\left(  M^{n}\times\mathbb{R}\right)
$ may be viewed as the Wall finiteness obstruction of $W\left(  M^{n}%
,N\right)  $ where $N$ is any sufficiently small $0$-neighborhood of infinity
in $M^{n}$, this obstruction $\sigma\left(  W\left(  M^{n},N\right)  \right)
$ lies in $\widetilde{K}_{0}\left(  \mathbb{Z[}\pi_{1}\left(  W\left(
N,m\right)  \right)  ]\right)  $. The retraction of $W\left(  N,m\right)  $
onto $M^{n}\times\left\{  m\right\}  $ obtained by projection shows that
$\pi_{1}(M^{n})$ is a retraction of $\pi_{1}\left(  W\left(  N,m\right)
\right)  $; so, by functorality, $\widetilde{K}_{0}\left(  \mathbb{Z[}\pi
_{1}(M^{n})]\right)  $ is a retraction of $\widetilde{K}_{0}\left(
\mathbb{Z[}\pi_{1}\left(  W\left(  N,m\right)  \right)  ]\right)  $. As a
consequence, the inclusion induced homomorphism $i_{\ast}:\widetilde{K}%
_{0}\left(  \mathbb{Z[}\pi_{1}(M^{n})]\right)  \rightarrow\widetilde{K}%
_{0}\left(  \mathbb{Z[}\pi_{1}\left(  W\left(  N,m\right)  \right)  ]\right)
$ is injective. By applying the Sum Theorem for Wall's finiteness obstruction
\cite[Ch.VI]{Si} to the homotopy equivalence $W\left(  N,m\right)  \simeq
K_{N}=M^{n}\cup_{\varphi}C^{\prime}$ utilized in the proof of Proposition
\ref{main-prop}, it is easy to see that $\sigma\left(  W\left(  M^{n}%
,N\right)  \right)  $ is precisely $i_{\ast}\left(  \sigma\left(
M^{n}\right)  \right)  $.

As a consequence of all of the above, we have a recipe for creating open
manifolds that are reasonably nice at infinity (inward tame with stable
fundamental group), but are not compactifiable via the addition of a manifold
boundary. Specifically: build a finite dimensional complex $K$ that is
finitely dominated but $\sigma\left(  K\right)  \neq0$; properly embed $K$ in
$\mathbb{R}^{n}$ and let $M^{n}$ be the interior of a proper regular
neighborhood of $K$; then $M^{n}\times\mathbb{R}$ satisfies Conditions 1) and
2) of Theorem \ref{Section: Proof}, but $\sigma_{\infty}(M^{n}\times
\mathbb{R)}$ is non-trivial and equal to $i_{\ast}\left(  \sigma\left(
K\right)  \right)  $.


\begin{thebibliography}{9999}                                                                                             %


\bibitem[Co]{Co}D.E. Cohen, \emph{Combinatorial Group Theory: a topological
approach}, London Mathematical Society, Student Texts 14, Cambridge University
Press, 1989.

\bibitem[Fr]{Fr}M.H. Freedman, \emph{The topology of four-dimensional
manifolds}, J. Differential Geom. 17 (1982), no. 3, 357--453.

\bibitem[Ge]{Ge}R. Geoghegan, \emph{Topological Methods in Group Theory},
book, writing in progress.

\bibitem[Gl]{Gl}J. Glimm, \emph{Two Cartesian products which are Euclidean
spaces}, Bull. Soc. Math. France 88 1960 131--135.

\bibitem[Gu]{Gu}C.R. Guilbault, \emph{Manifolds with non-stable fundamental
groups at infinity}, Geometry and Topology \textbf{4 }(2000), 537-579.

\bibitem[GuTi]{GuTi}C.R. Guilbault, \emph{Manifolds with nonstable fundamental
groups at infinity, II}, Geometry \& Topology, Volume 7 (2003), 255-286.

\bibitem[Hu]{Hu}S.T. Hu, \emph{Theory of retracts}, Wayne State University
Press, Detroit, 1965.

\bibitem[KS]{KS}R.C. Kirby and L.C. Siebenmann, \emph{Foundational essays on
topological manifolds, smoothings, and triangulations}, Annals of Mathematics
Studies, No. 88. Princeton University Press, Princeton, N.J., 1977. vii+355 pp.

\bibitem[Lu1]{Lu1}E. Luft, \emph{On contractible open topological manifolds},
Invent. Math. 4 1967 192--201.

\bibitem[Lu2]{Lu2}E. Luft, \emph{On contractible open} $3$\emph{-manifolds},
Aequationes Math. 34 (1987), no. 2-3, 231--239.

\bibitem[Ma]{Ma}W.S. Massey, \emph{Algebraic topology: an introduction},
Graduate Texts in Mathematics, Springer-Verlag, New York, 1987.

\bibitem[Mc]{Mc}D.R. McMillan, Jr. \emph{Cartesian products of contractible
open manifolds},.Bull. Amer. Math. Soc. 67 1961 510--514.

\bibitem[Si]{Si}L.C. Siebenmann, \emph{The obstruction to finding a boundary
for an open manifold of dimension greater than five}, Ph.D. thesis, Princeton
University, 1965.

\bibitem[St]{St}J. Stallings, \emph{The piecewise-linear structure of
Euclidean space}, Proc. Cambridge Philos. Soc. 58 1962, 481--488.

\bibitem[Wa]{Wa}C.T.C. Wall, \emph{Finiteness conditions for CW complexes},
Ann. Math. \textbf{8} (1965), 55-69.

\bibitem[Wh]{Wh}J.H.C. Whitehead, \emph{A certain open manifold whose group is
unity}, Quarterly J. Math., \textbf{6 }(1935), 268-279.
\end{thebibliography}
\end{document}